\documentclass[10pt,letterpaper]{amsart}
\usepackage{amssymb}
\usepackage{mathrsfs}
\usepackage{amsfonts}
\usepackage[active]{srcltx}

\makeatletter
\newtheorem{thm}{\textbf Theorem}[section]

\newtheorem{rem}{\textbf Remark}[section]

\newtheorem{prop}{\textbf Proposition}[section]

\newcommand{\md}{\mbox{d}}
\newcommand{\be}{\begin{equation}} 
\newcommand{\ee}{\end{equation}} 
\newcommand{\mr}{\mathbb{R}}

\newcommand{\mx}{\mbox}

\newcommand{\pt}{\partial}
\newcommand{\bes}{\begin{equation*}}
\newcommand{\ees}{\end{equation*}}

\newcommand{\al}{\alpha}

\newcommand{\La}{\Lambda}
\newcommand{\De}{\Delta}

\begin{document}
\begin{titlepage}
\title[Blowup criterion of incompressible viscoelastic flow]{Note On the blowup criterion of smooth solution to the incompressible viscoelastic flow}
\author{Baoquan Yuan}
\address{School of Mathematics and Informatics,
        Henan Polytechnic University, Henan, 454000, China.
        Courant Institute of Mathematical Sciences, New York University, New York, NY 10012 USA.}
\email{bqyuan@hpu.edu.cn}

\date{}
\end{titlepage}
\maketitle

\begin{abstract}
 {We study the blowup criterion of smooth solution to the Oldroyd
 models. Let $(u(t,x), F(t,x)$ be a smooth solution in $[0,T)$, it is
 shown that the solution $(u(t,x), F(t,x)$ does not appear breakdown until $t=T$ provided $\nabla u(t,x)\in L^1([0,T]; L^\infty(\mr^n))$, $n=2,3$.
}

\end{abstract}

 \vspace{.2in} {\bf  {Key words and phrases:}}{
Incompressible viscoelastic fluids, Oldroyd model, blowup criterion
of smooth solution.}


\section{{Introduction}}
 \setcounter{equation}{0}
In this paper, we consider the blowup criterion of smooth solution
to the incompressible Oldroy model in the two and three dimensional
space:
 \be\label{Od-B}
\begin{cases}
\pt_tu-\nu\De u+u\cdot\nabla u+\nabla p=\nabla\cdot(FF^t),\\
\pt_t F+u\cdot\nabla F=\nabla uF,\\
\mx{div} u=0,
\end{cases}
 \ee
for any $t>0$, $x\in \mr^n,$ $n=2,\ 3$, where $u(t,x)$ is the
velocity field, $p$ is the pressure, $\mu$ is the viscosity and $F$
the deformation tensor. We denote $(\nabla\cdot
F)_i=\pt_{x_j}F_{ij}$ for a matrix $F$. The Oldroy model
(\ref{Od-B}) describes an incompressible non-Newtonian fluid, which
bears the elastic property. For the details on this model see
\cite{Lin-Liu-Zhang}.

The local existence and uniqueness of the Oldroy model on entire
space $\mr^n$ or a periodic domain was established by Lin etc. in
\cite{Lin-Liu-Zhang}, where the global existence and uniqueness of
smooth solution with small initial data was also established see
also \cite{Lei-Liu-Zhou}. The wellposedness on a bounded smooth
domain with Dirichlet conditions was established by Lin and Zhang in
\cite{Lin-Zhang}.

We remark some properties of the deformation tensor. Let $x$ be the
Euler coordinate and $X$ the Lagrangian coordinate. For a given
velocity field $u(t,x)$ the flow map $x(t,X)$ is defined by the
following ordinary differential equation
 \bes
 \begin{cases}
\frac\md{\md t}x(t,X)=u(t,x(t,X)),\\
x(0,X)=X.
 \end{cases}
 \ees
The deformation tensor is $\tilde F(t,X)=\frac{\pt x}{\pt X}(t,X)$.
In the Eulerian coordinate, the corresponding deformation tensor is
define as $F(t,x(t,X))=\tilde F(t,X)$. Differentiating its both
sides with respect to $t$ by chain rule one obtain the second
equation of (\ref{Od-B}), which says that $\pt_tF_{ij}+u_k\cdot
\pt_{x_k}F_{ij}=\pt_{x_k}u_iF_{kj}$ for $i,\ j=1,2,\cdots, n$, in
the $(i,j)-$th entries, where we use the Einstein summation
convention that the repetition index denotes sum over $1$ to $n$.

If div $F(0,x)=0$, then from the second equation of Oldroy
(\ref{Od-B}) we have
 \be
 \pt_t(\nabla\cdot F^t)+u\cdot\nabla(\nabla\cdot F^t)=0.
 \ee
Therefore, $\nabla\cdot F^t=0$ for any $t> 0$.

Denote the $i$th column of $F$ as $F_{\cdot i}$, then $\nabla\cdot
(FF^t)=F_{\cdot i}\cdot\nabla F_{\cdot i}$ by the fact $\nabla\cdot
F^t=0$. So the system (\ref{Od-B}) can be rewritten in an equivalent
form
 \be\label{Od-B1}
\begin{cases}
\pt_tu-\nu\De u+u\cdot\nabla u+\nabla p=F_{\cdot i}\cdot\nabla F_{\cdot i},\\
\pt_t F_{\cdot k}+u\cdot\nabla F_{\cdot k}=F_{\cdot k}\cdot\nabla u,\ k=1,\cdots, n,\\
\mx{div} u=0,\ \mx{div} F=0.
\end{cases}
 \ee
In reference \cite{Lin-Liu-Zhang}, Lin, Liu and Zhang obtained the
local existence and uniqueness of smooth solution for smooth initial
data, and had a blowup criterion.

{\bf Theorem} (Lin, Liu and Zhang) For smooth initial data
$(u_0,F_0)\in H^2(\mr^n)$, there exists a positive time
$T=T(\|u_0\|_{H^2},\|F_0\|_{H^2})$ such that the system (\ref{Od-B})
possesses a unique smooth solution on $[0,T]$ with
 \bes
(u,F)\in L^\infty([0,T];H^2(\mr^n))\cap L^2([0,T];H^3(\mr^n)).
 \ees
Moreover, if $T^*$ is the maximal time of existence, then
 \bes
\int^{T^*}_0\|\nabla u\|^2_{H^2}\md s=+\infty.
 \ees
In reference \cite{Hu-Hynd}, Hu and Hynd study the blowup criterion
for the ideal viscoelastic flow, which is the Oldroy system
(\ref{Od-B}) in the case of $\mu=0$. They showed an Beale-Kato-Majda
\cite {Beale-Kato-Majda} type blowup criterion that the smooth
solution to the Oldroy flow do not develop singularity for $t\le T$
provided that
 \bes
\int^T_0\|\nabla\times u\|_{L^\infty(\mr^3)}\md
s+\sum_{k=1}^3\int^T_0\|\nabla\times F_{\cdot
k}\|_{L^\infty(\mr^3)}\md s<+\infty.
 \ees
From the modeling of Oldroy system we know that the deformation
tensor can be determined by the velocity $u$ of the flow. Therefore
we consider the blowup criterion of smooth solution by means of only
$\|\nabla u\|_{\infty}$. In fact, Zhao, Guo and Huang
\cite{Zhao-Guo-Huang} constructed a set of finite time blowup
solution in two dimension case:
 \bes
 \begin{split}
&u(t,x)=\bigg(\frac{x_1 f_0}{1-\frac{\alpha+\beta}{\al-\beta}f_0t},
\frac{x_2 f_0}{1-\frac{\alpha+\beta}{\al-\beta}f_0t}\bigg)^t, \
p(t,x)=\frac{(\al x_1^2-\beta
x_2^2)f_0^2}{(\beta-\al)(1-\frac{\al+\beta}{\al-\beta}f_0t)^2},\\
 &F(t,x)=\mx{diag}\bigg(\Big|1-\frac{\al+\beta}{\al-\beta}f_0t\Big|^{\frac{\beta-\al}{\al+\beta}},
  \Big|1-\frac{\al+\beta}{\al-\beta}f_0t\Big|^{\frac{\beta+\al}{\al-\beta}}\bigg).
 \end{split}
 \ees
If $\frac{\al+\beta}{\al-\beta}f_0>0$, $\al+\beta\neq 0$ and
$\al-\beta\neq 0$, then the above solution will blow up at time
$T^*=\frac{\al-\beta}{(\al+\beta)f_0}$. We see that
 \bes
\int^{T^*}_0\|\nabla u(t)\|_{\infty}\md t=+\infty.
 \ees
There are other types of blowup criteria of smooth solutions to the
Oldroy models, for example \cite{Lei-Masmoudi-Zhou,
Chemin-Masmoudi}. To this end, we state our main results.

\begin{thm}\label{Thm1}
Let $u_0\in H^2(\mr^n)$ and $F_0\in H^2(\mr^n)$ with $\nabla\cdot
u_0=\nabla\cdot F_{\cdot k,0}=0$ for $k=1,\cdots, n$. Assume the
pair $(u,F)\in L^\infty([0,T];H^2(\mr^n))\cap L^2([0,T];H^3(\mr^n))$
is a smooth solution to the Oldroy system (\ref{Od-B1}). Then the
smooth solution do not appear breakdown until $T^*> T$ provided that
 \be
\int^{T^*}_0\|\nabla u(t)\|_{\infty}\md t<+\infty.
 \ee
\end{thm}

\begin{rem}
For the local smooth solution $(u,F)\in
L^\infty([0,T];H^2(\mr^n))\cap L^2([0,T];H^3(\mr^n))$, if $T^*$ is
its maximum existence time, then $\int^{T^*}_0\|\nabla
u(t)\|_{\infty}\md t=+\infty$.
\end{rem}

In the second section we will prove the Theorem \ref{Thm1} for the
case $n=2$, which can be done by energy estimates. The $L^2$ and
$H^1$ energy estimates are the same for the case $n=2$ and $n=3$. In
the $H^2$ energy estimate, we use the Sobolev interpolation
inequality $\|\nabla F\|^2_4\le C\|\nabla F\|_2\|\De F\|_2$. In case
$n=3$, however, the inequality is $\|\nabla F\|^2_4\le C\|\nabla
F\|_2^{\frac12}\|\De F\|_2^{\frac32}$ which does not match the $H^2$
energy estimate, because it will result in the appearance of the
term $\|\De F\|^3_2$ that the power is higher that the left hand
side. We obtain the $H^2$ energy estimate of $u$ by virtue of the
momentum equation, combining the $H^2$ estimate of $u$ and $F$ again
with the estimate of $\|\nabla F\|_{L^6}$ we grasp the $H^2$ energy
estimate of $u$ and $F$ finally. The section three will devote to
the proof of the case $n=3$.

In this paper $C$ denote a harmless constant which may be dependent
on dimension $n$, the norm of initial data, the viscosity $\mu$, but
not dependent on the estimated quantity. We denote the $L^p$ norm of
a function $f$ by $\|f\|_p$ or $\|f\|_{L^p}$. We denote the
derivative with respect to $x_i$ by $\pt_i$ or $\pt_{x_i}$. We also
use $f_t$ to denote the derivative of $f$ with respect to $t$.


\section {Proof of the case $n=2$}
\setcounter{equation}{0}

(1) $L^2$-energy estimate and $L^p$ estimate of the deformation
tensor $F$

The $L^2$-energy estimate can be easily obtained by the standard
$L^2$ inner product process.
 \bes
\frac12\frac{\md}{\md t}(\|u\|_2^2+\|F_{\cdot k}\|^2_2)+\mu\|\nabla
u\|^2_2=(F_{\cdot i}\cdot\nabla F_{\cdot i},u)+(F_{\cdot
k}\cdot\nabla u,F_{\cdot i})=0.
 \ees
So we have
 \be\label{Energy1}
\|u\|^2_2+\|F\|^2_2+2\mu\int^t_0\|\nabla u\|_2^2\md
s=\|u_0\|_2^2+\|F(0)\|_2^2.
 \ee

Multiplying both sides of the second equation of (\ref{Od-B1}) by
$p|F_{\cdot k}|^{p-2}F_{\cdot k}$ for $2\le p<\infty$ and
integrating both sides on $\mr^n$ it follows that
 \be\label{Lp1}
\frac \md{\md t}\|F_{\cdot k}\|_p^p\le p\|\nabla
u\|_\infty\|F\|_p^p.
 \ee
Summing up the estimate (\ref{Lp1}) with respect to k one has
 \be\label{Lpestimate}
\|F\|_p\le \|F_0\|_p\exp\bigg\{C(n)\int^t_0\|\nabla u(s)\|_\infty\md
s\bigg\}.
 \ee
 Let $p\rightarrow \infty$, we have
 \be\label{Linfinityestimate}
 \|F\|_\infty\le
\|F_0\|_\infty\exp\bigg\{C(n)\int^t_0\|\nabla u(s)\|_\infty\md
s\bigg\}.
 \ee

 (2) $\dot H^1$-energy estimate

We differentiate the equations (\ref{Od-B1}) with respect to $x_i$,
then multiply the resulting equations by $\pt_i u$ and $\pt_i
F_{\cdot j}$ for $i = 1,\ 2$, integrate with respect to $x$ and sum
them up. It follows that

 \bes
 \begin{split}
&\frac12\frac
\md{\md t}(\|\pt_iu\|^2_2+\|\pt_iF\|^2_2)+\mu\|\pt_i\nabla u\|^2_2\le\\
& |(\pt_iu\cdot\nabla u,\pt_iu)|+|(\pt_iF_{\cdot k}\cdot\nabla
F_{\cdot k},\pt_iu)| +|(\pt_iu\cdot\nabla F_{\cdot j},\pt_i F_{\cdot
j})|+|(\pt_iF_{\cdot j}\cdot\nabla u, \pt_iF_{\cdot j})|,
\end{split}
 \ees
where use has been made of the facts
 \bes
 \begin{split}
&(u\cdot\nabla\pt_iu,\pt_iu)=(u\cdot\pt_iF_{\cdot j},\pt_iF_{\cdot
j})=(\nabla\pt_i p,\pt_i u)=0,\\
&(F_{\cdot k}\cdot\nabla\pt_i F_{\cdot k},\pt_i u)+(F_{\cdot
j}\cdot\nabla\pt_i u,\pt_i F_{\cdot j})=0.
 \end{split}
 \ees

 Noting that
 \bes
|(\pt_iu\cdot\nabla u,\pt_iu)|\le \|\nabla u\|_{\infty}\|\nabla
u\|^2_2,
 \ees
 \bes
|(\pt_iF_{\cdot k}\cdot\nabla F_{\cdot k},\pt_iu)|,\
|(\pt_iu\cdot\nabla F_{\cdot j},\pt_iF_{\cdot j})|,\ |(\pt_iF_{\cdot
j}\cdot\nabla u,\pt_iF_{\cdot j})|\le\|\nabla u\|_{\infty}\|\nabla
F\|^2_2.
 \ees
So
 \bes
\frac12\frac \md{\md t}(\nabla u\|^2_2+\|\nabla
F\|^2_2)+\mu\|D^2u\|^2_2\le C\|\nabla u\|_\infty(\|\nabla
u\|^2_2+\|\nabla F\|^2_2).
 \ees
Gronwall's inequality implies
 \be\label{Energy2}
\|\nabla u\|^2_2+\|\nabla F\|^2_2+2\mu\int^t_0\|D^2u\|^2_2\md s\le
(\|\nabla u_0\|^2_2+\|\nabla F(0)\|^2_2)\exp\Bigg\{\int^t_0C\|\nabla
u(s)\|_\infty\md s\Bigg\}.
 \ee

(3) $\dot H^2$-energy estimate

Applying operator $\De$ on both sides of (\ref{Od-B1}), we have
 \be\label{E2L}
 \begin{cases}
\pt_t\De u-\mu\De^2u+\De u\cdot\nabla u+u\cdot\nabla\De
u+2\pt_iu\cdot\nabla\pt_i u+\nabla\De p=\De F_{\cdot k}\cdot\nabla
F_{\cdot k}+ F_{\cdot k}\nabla\De F_{\cdot k}+2\pt_i
F_{\cdot k}\cdot\nabla\pt_i F_{\cdot k}\\
\pt_t\De F_{\cdot k}+\De u\cdot\nabla F_{\cdot k}+u\cdot\nabla\De
F_{\cdot k}+2\pt_i u\cdot\nabla\pt_i F_{\cdot k}=\De F_{\cdot
k}\cdot\nabla u+F_{\cdot k}\cdot\nabla\De u+2\pt_i F_{\cdot
k}\cdot\nabla\pt_iu.
 \end{cases}
 \ee
Taking the $L^2$ inner of equation (\ref{E2L}) with $\De u$ and $\De
F_{\cdot k}$ and summing them up, one can obtain that
 \be\label{E3}
 \begin{split}
&\frac12\frac \md{\md t}(\|\De u\|^2_2+\|\De F\|^2_2)+\mu\|\De\nabla
u\|^2_2 \\ &\le |(\De u\cdot\nabla u,\De u)|+2|(\pt_i
u\cdot\nabla\pt_i
u,\De u)|+|(\De F_{\cdot k}\cdot\nabla F_{\cdot k},\De u)|\\
&+2|(\pt_i F_{\cdot k}\cdot\nabla \pt_i F_{\cdot k},\De u)|+|(\De
u\cdot\nabla F_{\cdot k}, \De F_{\cdot k})|+2|(\pt_i
u\cdot\nabla\pt_i F_{\cdot k}, \De F_{\cdot k})|\\ &+|(\De F_{\cdot
k}\cdot\nabla u, \De F_{\cdot k})|+2|(\pt_i F_{\cdot
k}\cdot\nabla\pt_i u, \De F_{\cdot k})|.
 \end{split}
 \ee
Here use has been made of the the facts that
 \bes
 \begin{split}
(u\cdot\nabla\De u,\De u)=0, \ (u\cdot\nabla\De F_{\cdot k},\De
F_{\cdot k})=0,\\
(F_{\cdot k}\cdot\nabla\De F_{\cdot k},\De u)+(F_{\cdot
k}\cdot\nabla\De u,\De F_{\cdot k})=0.
 \end{split}
 \ees
Noting that
 \bes
|(\De u\cdot\nabla u, \De u)|,\  |(\pt_iu\cdot\nabla\pt_i u,\De
u)|\le C\|\nabla u\|_\infty\|\De u\|^2_2,
 \ees
 \bes
|(\De F_{\cdot k}\cdot\nabla u, \De F_{\cdot k})|,\ |(\pt_i
u\cdot\nabla\pt_i F_{\cdot k}, \De F_{\cdot k})|\le C\|\nabla
u\|_\infty\|\De F\|^2_2.
 \ees
 \bes
 \begin{split}
&|(\De F_{\cdot k}\cdot\nabla F_{\cdot k},\De u)+2(\pt_i F_{\cdot
k}\cdot\nabla\pt_i F_{\cdot k},\De u)|\\ &=|-(\pt_i F_{\cdot
k}\cdot\nabla F_{\cdot k},\pt_i\De u)-(\pt_i F_{\cdot
k}\cdot\nabla\De u,\pt_i F_{\cdot k})|\le C\|\nabla \De u\|_2\|\nabla F\|^2_4\\
&\le \frac \mu4\|\nabla \De u\|_2^2+ C\|\nabla F\|_2^2\|\De F\|^2_2,
\end{split}
 \ees
where we have used the Sobolev interpolation inequality
 \bes
\|\nabla F\|^2_4\le C\|\nabla F\|_2\|\De F\|_2.
 \ees
Arguing similarly as the above, one has
 \bes
|(\De u\cdot\nabla F_{\cdot k},\De F_{\cdot k})| =|(\pt_i\De
u\cdot\nabla F_{\cdot k},\pt_i F_{\cdot k})|\le \frac \mu
8\|D^3u\|_2^2+ C\|\nabla F\|_2^2\|\De F\|^2_2,
 \ees
 \bes
 \begin{split}
&|(\pt_i F_{\cdot k}\cdot\nabla\pt_i u,\De F_{\cdot
k})|=|(\pt_i\pt_j F_{\cdot k}\cdot\nabla\pt_j F_{\cdot k},\pt_i
u)|+|(\pt_i F_{\cdot k}\cdot\nabla\pt_i\pt_j u,\pt_j F_{\cdot k})|\\
&\le C\|\nabla u\|_\infty\|\De F\|_2^2+\frac \mu 8\|\nabla \De
u\|_2^2+ C\|\nabla F\|_2^2\|\De F\|^2_2.
 \end{split}
 \ees
Inserting the above estimates into estimate (\ref{E3}), it can be
derived that
 \bes
\frac12\frac \md{\md t}(\|\De u\|^2_2+\|\De F\|^2_2)+\frac\mu
2\|\De\nabla u\|^2_2\le C\|\nabla u\|^2_\infty(\|\De u\|^2_2+\|\De
F\|^2_2)+C\|\nabla F\|_2^2\|\De F\|^2_2.
 \ees
Gronwall's inequality implies that
 \be\label{Energy3}
\|\De u\|^2_2+\|\De F\|^2_2+\mu \int^t_0\|\De\nabla u(s)\|^2_2\md
s\le (\|\De u_0\|^2_2+\|\De F(0)\|^2_2)\exp\Bigg\{ C\exp \int^t_0
Ct\|\nabla u\|_\infty\md s \Bigg\}.
 \ee

(4) Higher derivative estimates.

Next we derive the higher derivative estimate of $u$ and $F$. For
this purpose we need the following commutator estimate.

\begin{prop} {\rm (Kato and Ponce \cite{kato-Ponce}, \cite{Majda-Bertozzi})}
Let $1<p<\infty$ and $0<s$. Assume that $f,\ g\in W^{s,p}$, then
there exists a abstract constant $C$ such that
 \be\label{CommutatorEstimate1}
\|[J^s,f]g\|_p\le C(\|\nabla
f\|_{p_1}\|g\|_{W^{s-1,p_2}}+\|f\|_{W^{s,p_3}}\|g\|_{p_4})
 \ee
 \be\label{CommutatorEstimate}
\|[\La^s,f]g\|_p\le C(\|\nabla f\|_{p_1}\|\La^{s-1}g\|_{p_2}+\|\La
^sf\|_{p_3}\|g\|_{p_4})
 \ee

 with $1<p_2, p_3<\infty$ such that
 $$\frac1p=\frac1{p_1}+\frac1{p_2}=\frac1{p_3}+\frac1{p_4},$$ where
 $[\La^s,f]g=\La^s(fg)-f\La^sg$ and $\La^s=(-\De)^{\frac12}$, $J=(1-\De)^{1/2}$.
 \end{prop}

Applying $\La^s$ on both sides of (\ref{Od-B1}) and taking the inner
product with $\La^s u$ and $\La^s F$, it can be derived that
 \be\label{Higher1}
 \begin{split}
&\frac12\frac\md {\md t}(\|\La^s u\|_2^2+\|\La^s F_{\cdot
k}\|_2^2)+\mu \|\La^{s+1}u\|_2^2\le \\&|(\La^s(u\cdot \nabla
u)-u\cdot\nabla\La^s u,\La^s u)|+|(\La^s(F_{\cdot k}\cdot \nabla
F_{\cdot k})-F_{\cdot k}\cdot\nabla\La^s F_{\cdot k},\La^s u)|+ \\
&|(\La^s(u\cdot \nabla F_{\cdot k})-u\cdot\nabla\La^s F_{\cdot
k},\La^s F_{\cdot k})|+|(\La^s(F_{\cdot k}\cdot \nabla u)-F_{\cdot
k}\cdot\nabla\La^s u,\La^s F_{\cdot k})|,
 \end{split}
 \ee
where we have used the facts
 \bes
(F_{\cdot k}\cdot\nabla\La^s F_{\cdot k}, \La^s u)+(F_{\cdot
k}\cdot\nabla\La^s u, \La^s F_{\cdot k})=0,
 \ees
 \bes
(u\cdot\nabla\La^s F_{\cdot k}, \La^s F_{\cdot
k})=(u\cdot\nabla\La^s u, \La^s u)=0.
 \ees
The commutator estimate (\ref{CommutatorEstimate}) implies that
 \bes
\|\La^s(u\cdot\nabla u)-u\cdot\nabla\La^s u\|_2\le \|\nabla
u\|_\infty\|\La^s u\|_2,
 \ees
 \bes
\|\La^s(F_{\cdot k}\cdot\nabla F_{\cdot k})-F_{\cdot
k}\cdot\nabla\La^s F_{\cdot k}\|_2\le \|\nabla F\|_{\infty}\|\La^s
F\|_2\le \|\nabla F\|_{H^{s-1}}\|\La^s F\|_2,
 \ees
 \bes
\|\La^s(u\cdot\nabla F_{\cdot k})-u\nabla\La^s F_{\cdot k}\|\le
\|\nabla u\|_\infty\|\La^s F\|_2+\|F\|_\infty\|\La^{s+1}u\|_2,
 \ees
 \bes
 \|\La^s(F_{\cdot k}\cdot\nabla u)-F_{\cdot k}\nabla\La^s u\|\le
\|\nabla u\|_\infty\|\La^s F\|_2+\|F\|_\infty\|\La^{s+1}u\|_2,
 \ees
where the Sobolev embedding $H^{s-1}(\mr^n)\hookrightarrow
L^{\infty}(\mr^n)$ for $s>1+\frac n2$ is applied.

Inserting the above estimates into estimate (\ref{Higher1}), it
follows
 \be\label{Higher2}
 \begin{split}
&\frac12\frac\md {\md t}(\|\La^s u\|_2^2+\|\La^s F_{\cdot
k}\|_2^2)+\frac\mu2 \|\La^{s+1}u\|_2^2\le \\&C(\|\nabla
u\|_\infty+\|\nabla F\|_2+\|\La^s u\|_2+\|F\|_\infty^2)(\|\La^s
u\|_2^2+\|\La^s F\|_2^2),
 \end{split}
 \ee
where we have used the fact
 \bes
\|\nabla F\|_{H^{s-1}}\|\La^s F\|_2\|\La^s u\|_2\le \|\nabla
F\|_2(\|\La^s F\|^2_2+\|\La^s u\|_2^2)+\|\La^s F\|_2^2\|\La^s u\|_2.
 \ees
So, for $s\ge 3$, applying Gronwall's inequality to (\ref{Higher2}),
by induction for $u$'s estimate, we obtain the higher derivative
estimate:
 \bes
 \begin{split}
&\|\La^s u\|_2^2+\|\La ^sF\|_2^2+\mu\int^t_0\|\La^{s+1}u\|_2^2\md
s\le
\\ &(\|u_0\|_{H^s}^2+\|F(0)\|_{H^s}^2)\exp\bigg\{\int^t_0C(\|\nabla
u\|_\infty+\|\nabla F\|_2+\|\La^s u\|_2+\|F\|_\infty^2)\md s\bigg\}.
 \end{split}
 \ees
Therefore, we complete the proof of the case $n=2$.


\section {Proof of the case $n=3$}
\setcounter{equation}{0}

In the three dimensional case the $L^2$ and $H^1$ energy estimates
are the same as the case of dimension two. To estimate the $H^2$
energy estimate we need the following estimates.

Multiplying the first equation of (\ref{Od-B1}) by $u_t$ and
integrating both sides over $\mr^3$ with respect to $x$, and noting
div $u=0$, it follows
 \bes
 \begin{split}
\frac\mu 2\frac\md{\md t}\|\nabla u\|_2^2+\|u_t\|_2^2
&\le|(u\cdot\nabla u, u_t)|+|(F_{\cdot k}\cdot\nabla F_{\cdot
k},u_t)|\\ &\le \frac12\|u_t\|_2^2+C\|u\|_\infty^2\|\nabla
u\|_2^2+C\|\nabla F\|_2^2\|F\|_\infty^2.
 \end{split}
 \ees
Integrating both sides with respect to $t$ it yields
 \be\label{Testimate1}
\mu\|\nabla u\|_2^2+\int^t_0\|u_t\|^2_2\md s\le\mu\|\nabla
u_0\|_2^2+\sup_{0<s<t}\|\nabla u\|_2^2\int^t_0\|u\|^2_{H^2}\md
s+\int^t_0\|F\|^2_\infty\|\nabla F\|_2^2\md s
 \ee
 where the Sobolev embedding $H^2(\mr^3)\hookrightarrow
 L^\infty(\mr^3)$ has been used.

Differentiating the first equation of (\ref{Od-B1}) with respect to
$t$, we arrive at
 \be\label{tteqs}
u_{tt}-\mu\De u_t+u_t\cdot\nabla u+u\cdot\nabla u_t+\nabla
p_t=F_{\cdot kt}\cdot\nabla F_{\cdot k}+F_{\cdot k}\cdot\nabla
F_{\cdot kt}.
 \ee
Taking $L^2$ inner product of the equation (\ref{tteqs}) with
respect to $u_t$, it can be similarly derived that
 \bes
 \begin{split}
\frac12\frac\md{\md t}\|u_t\|_2^2+\mu\|\nabla u_t\|_2^2&\le \|\nabla
u\|_\infty\|u_t\|_2^2+2\|F\|_\infty\|\nabla u_t\|_2\|F_t\|_2\\ &\le
\frac\mu 2\|\nabla u_t\|_2^2+\|\nabla
u\|_\infty\|u_t\|^2_2+C\|F\|^2_\infty\|F_t\|_2^2.
\end{split}
 \ees
Applying the Gronwall' inequality, it yields
 \be\label{ttestimate}
\|u_t\|_2^2+\mu\int^t_0\|\nabla u_t\|_2^2\md s\le
(\|u_t(0)\|_2^2+C\int^t_0\|F\|_\infty^2\|F_t\|_2^2\md
s)\exp\bigg\{\int^t_0\|\nabla u\|_\infty\md s\bigg\}.
 \ee
It need still to estimate $\|F_t\|_2^2$. From the second equation of
(\ref{Od-B1}) it can be derived that
 \bes
 \begin{split}
\|F_t\|_2^2&\le \|F_t\|_2\|u\|_\infty\|\nabla
F\|_2+\|F_t\|_2\|F\|_\infty\|\nabla u\|_2\\ &\le
\frac12\|F_t\|_2^2+C\|u\|_\infty^2\|\nabla
F\|^2_2+C\|F\|^2_\infty\|\nabla u\|_2^2.
 \end{split}
 \ees
So we arrive at
 \bes
\|F_t\|_2^2\le C\|u\|_\infty^2\|\nabla
F\|^2_2+C\|F\|^2_\infty\|\nabla u\|_2^2.
 \ees
Inserting it to the estimate (\ref{ttestimate}) we obtain the
estimate of $\|u_t\|_2$:
 \be\label{ttestimate1}
\|u_t\|_2^2+\mu\int^t_0\|\nabla u_t\|_2^2\md s\le C(t)<\infty,
 \ee
where $C(t)$ is explicit increasing function of $t$ dependent on
$\int^t_0\|\nabla u\|_\infty\md s$. From the first equation of
(\ref{Od-B1}), $\nabla p$ can be solved by Riesz transformation
$R=(R_1,R_2,R_3)^t$, with $R_j=-i\pt_{x_j}(-\De)^{-\frac12}$ being
the $j$th Riesz transformation.
 \bes
\nabla p=RR\cdot(u\cdot \nabla u)-RR\cdot(F_{\cdot k}\cdot\nabla
F_{\cdot k}).
 \ees
In virtue of the boundedness of Riesz operator R in $L^p$ space for
$1<p<\infty$, we obtain that
 \bes
\|\nabla p\|_2\le C\|\nabla u\|_2\|u\|_\infty+C\|\nabla
F\|_2\|F\|_\infty.
 \ees
For details about Riesz transformation see \cite{Miao,S-W}.

Thus from the first equation of (\ref{Od-B1}) we have
 \bes
 \begin{split}
\mu\|\De u\|_2&\le\|u_t\|_2+\|u\cdot \nabla u\|_2+\|\nabla
p\|_2+\|F_{\cdot k}\cdot\nabla F_{\cdot k}\|_2\\ &\le
\|u_t\|_2+\frac\mu2\|\De u\|_2+C\|u\|_2\|\nabla
u\|_2^4+C\|F\|_\infty\|\nabla F\|_2,
 \end{split}
 \ees
 where the interpolation inequality $\|u\|_\infty\le C\|u\|_2^{\frac14}\|\De
 u\|_2^{\frac34}$ has been used.
So we derive
 \be\label{uH2}
\|\De u\|_2\le C(\|u_t\|_2+\|u\|_2\|\nabla
u\|_2^4+\|F\|_\infty\|\nabla F\|_2).
 \ee

Next we derive the estimate of $\|\De F\|_2$. Applying $\De$ on the
both sides of equation (\ref{Od-B1}) and taking the $L^2$ inner
product with $\De u$ and $\De F_{\cdot k}$ respectively, we have
 \be\label{F1}
\frac12\frac\md{\md t}\|\De u\|_2^2+\mu\|\De\nabla u\|_2^2\le |(\De
(u\cdot\nabla u)-u\cdot\nabla\De u,\De u)|+|(\De (F_{\cdot
k}\cdot\nabla F_{\cdot k})-F_{\cdot k}\cdot\nabla\De F_{\cdot k},\De
u)|,
 \ee
 \be\label{F2}
\frac12\frac\md{\md t}\|\De F_{\cdot k}\|_2^2\le |(\De (u\cdot\nabla
F_{\cdot k})-u\cdot\nabla\De F_{\cdot k},\De F_{\cdot k})|+|(\De
(F_{\cdot k}\cdot\nabla u)-F_{\cdot k}\cdot\nabla\De u,\De F_{\cdot
k})|,
 \ee
where use has been of the facts
 \bes
(u\cdot\nabla\De u,\De u)=(u\cdot\nabla\De F_{\cdot k},\De F_{\cdot
k})=0,
 \ees
 \bes
(F_{\cdot k}\cdot\nabla\De F_{\cdot k},\De u)+(F_{\cdot
k}\cdot\nabla\De u,\De F_{\cdot k})=0.
 \ees
Next we estimate the right hand sides. By the communicator estimate
(\ref{CommutatorEstimate}) one has
 \bes
|(\De (u\cdot\nabla u)-u\cdot\nabla\De u,\De u)|\le \|\De u\|_2\|\De
(u\cdot\nabla u)-u\cdot\nabla\De u\|_2\le \|\nabla u\|_\infty\|\De
u\|_2^2,
 \ees
 \bes
 \begin{split}
&|(\De (u\cdot\nabla F_{\cdot k})-u\cdot\nabla\De F_{\cdot k},\De
F_{\cdot k})|\le \|\De F\|_2(\|\nabla u\|_\infty\|\De
F\|_2+\|F\|_\infty\|\nabla\De u\|_2)\\ &\le \|\nabla u\|_\infty\|\De
F\|_2^2+C\|F\|_\infty^2\|\De F\|_2^2+ \frac\mu{8}\|\nabla\De
u\|_2^2,
 \end{split}
 \ees
 \bes
|(\De (F_{\cdot k}\cdot\nabla u)-F_{\cdot k}\cdot\nabla\De u,\De
F_{\cdot k})|\le \|\nabla u\|_\infty\|\De
F\|_2^2+C\|F\|_\infty^2\|\De F\|_2^2+ \frac\mu{8}\|\nabla\De
u\|_2^2.
 \ees
For the second term on the right hand side of (\ref{F1}) we estimate
as follows
 \bes
|(\De (F_{\cdot k}\cdot\nabla F_{\cdot k})-F_{\cdot k}\cdot\nabla\De
F_{\cdot k},\De u)|\le \|\De u\|_6\|\De (F_{\cdot k}\cdot\nabla
F_{\cdot k})-F_{\cdot k}\cdot\nabla\De F_{\cdot k}\|_{6/5},
 \ees
and
 \bes
\|\De (F_{\cdot k}\cdot\nabla F_{\cdot k})-F_{\cdot k}\cdot\nabla\De
F_{\cdot k}\|_{6/5}\le \|\nabla F\|_6\|\De F\|_{3/2}\le \|\nabla
F\|_6\|\nabla F\|_2^{\frac12}\|\De F\|_2^{\frac12}.
 \ees
So one has the estimate
 \bes
|(\De (F_{\cdot k}\cdot\nabla F_{\cdot k})-F_{\cdot k}\cdot\nabla\De
F_{\cdot k},\De u)|\le \frac\mu 4\|\nabla\De u\|_2^2+C\|\nabla
F\|_6^4+C\|\nabla F\|_2^2\|\De F\|_2^2.
 \ees
Summing up (\ref{F1}) and (\ref{F2}), and inserting the above
estimates into the summation, we arrive at
 \be\label{F3}
\frac\md{\md t}(\|\De u\|_2^2+\|\De F_{\cdot
k}\|_2^2)+\mu\|\nabla\De u\|_2^2\le C(\|\nabla
u\|_\infty+\|F\|^2_\infty+\|\nabla F\|_2^2)(\|\De u\|_2^2+\|\De
F\|_2^2)+C\|\nabla F\|^4_6.
 \ee
We still have to estimate $\|\nabla F\|_6$. Differentiating the
second equation of (\ref{Od-B1}) with respect to $x_i$, one has
 \bes
\pt_t\pt_iF_{\cdot k}+\pt_iu\cdot\nabla F_{\cdot
k}+u\cdot\nabla\pt_i F_{\cdot k}=\pt_i F_{\cdot k}\cdot\nabla
u+F_{\cdot k}\cdot\nabla\pt_i u.
 \ees
Multiplying both sides of the above equation by $6|\pt_iF_{\cdot
k}|^4\pt_iF_{\cdot k}$, and integrating both sides with respect to
$x$ over $\mr^3$, it can be derived that
 \be\label{L6}
\frac\md{\md t}\|\nabla F\|_6^4\le C\|\nabla u\|_\infty\|\nabla
F\|_6^4+ C\|F\|_\infty\|\De u\|_6\|\nabla F\|_6^3.
 \ee
Next we have to derive an estimate of $\|\De u\|_6$. Using an
argument similar to deriving the $L^2$ estimate $\|\De u\|_2$ in
(\ref{uH2}) we have
 \be \label{u3}
 \begin{split}
\mu\|\De u\|_6 &\le \|\pt_tu\|_6+\|u\|_\infty\|\nabla
u\|_6+C\|F\|_\infty\|\nabla F\|_6\\ &\le \|\pt_t\nabla
u\|_2+C\|u\|_\infty\|\De u\|_2+C\|F\|_\infty\|\nabla F\|_6.
 \end{split}
 \ee
Inserting estimates (\ref{u3}) to (\ref{L6}) one has
 \be \label{L6(1)}
 \begin{split}
\frac\md {\md t}\|\nabla F\|_6^4 &\le C(\|\nabla
u\|_{\infty}+\|F\|^2_\infty+\|\pt_t \nabla u\|_2^2)\|\nabla
F\|_6^4+C\|F\|_\infty\|u\|_\infty\|\De u\|_2\|\nabla F\|_6^3+C\\
&\le C(\|\nabla u\|_{\infty}+\|F\|^2_\infty+\|\pt_t \nabla
u\|_2^2+1)\|\nabla F\|_6^4+C\|F\|_\infty^4\|u\|_2 \|\De u\|_2^7+C.
 \end{split}
 \ee
Combining the estimates (\ref{F3}) and (\ref{L6(1)}) we arrive at
 \bes
 \begin{split}
&\frac\md{\md t}(\|\De u\|_2^2+\|\De F_{\cdot
k}\|_2^2+\|F\|_6^4)+\mu\|\nabla\De u\|_2^2\le \\ &C(\|\nabla
u\|_\infty+\|\nabla F\|_2^2+\|F\|^2_\infty+\|\pt_t \nabla
u\|_2^2+1)(\|\De u\|_2^2+\|\De
F\|_2^2+\|F\|_6^4)+C\|F\|_\infty^4\|u\|_2 \|\De u\|_2^7+C.
 \end{split}
 \ees
Gronwall's inequality implies the $H^2$ estimates:
 \be\label{H2}
 \begin{split}
&\|\De u\|_2^2+\|\De F_{\cdot
k}\|_2^2+\|F\|_6^4+\mu\int^t_0\|\nabla\De u\|_2^2\md s\le
\exp\bigg\{C(t)\int^t_0(\|\nabla u\|_\infty+\|\pt_t\nabla
u\|_2^2)\md s\bigg\}\times \\ &\bigg(\|\De u(0)\|_2^2+\|\De F_{\cdot
k}(0)\|_2^2+\|F(0)\|_6^4+C\int^t_0(\|F\|_\infty^4\|u\|_2 \|\De
u\|_2^7+1)\md s\bigg)<\infty.
 \end{split}
 \ee
Based on the $H^2$ energy estimate the higher energy estimate can be
obtained by bootstrap method as we did in section two. Thus the
proof of the case $n=3$ is completed.

 \vspace{0.4cm}

\textbf{Acknowledgements} This work was done when the author was
visiting the Courant Institute of Mathematical Sciences at New York
University. The author would like to thank Professor Fanghua Lin for
stimulating discussion on this topic. The research of B Yuan was
partially supported by the National Natural Science Foundation of
China (No. 11071057), Innovation Scientists and Technicians Troop
Construction Projects of Henan Province (No. 104100510015), Program
for Science\&Technology Innovation Talents in Universities of Henan
Province (No. 2009 HASTIT007) and Doctor Fund of Henan Polytechnic
University (No. B2008-62).

\bibliographystyle{amsplain}

\end{document}